\documentclass[runningheads]{svjour2}
\smartqed  
\usepackage{graphicx}
\usepackage{latexsym,amssymb,amscd,amsmath}

\def \qed {\hfill \rule{1ex}{1ex}}
\def \U {{\bf u}}
\def \x {{\bf x}}

\def \H {{\bf H}}
\def \G {{\bf g}}

\def \V {{\bf v}}
\def \bt {{\widetilde{\hbox{\rm \bf b}}}}
\def \W {{\bf w}}
\def \N {{\mathbf{n}}}

\def \Ut {{\tilde {\mathbf{u}}}}
\def \Deltat {{\mathbf{\widetilde{\Delta}}}}
\def \F {{\bf f }}
\def \P {{\bf P}}
\def \R {{\mathbf{R}}}

\def \L {{\bf L}}

\def \b {{\hbox{\rm  b}}}

\def \Nabla {{\hbox{\boldmath $\nabla$ \unboldmath \!\!}}}
\def \Del {{\hbox{\boldmath $\delta$ \unboldmath \!\!}}}
\def \psig {{\hbox{\boldmath $\psi$ \unboldmath \!\!}}}
\def \alphag {{\hbox{\boldmath $\alpha$ \unboldmath \!\!}}}
\let \eps=\varepsilon

\newcommand{\Eta}[2]{\hbox{\boldmath $ \eta^{\hbox{\unboldmath $\scriptstyle#1$}}
_{\hbox{\unboldmath $\scriptstyle #2$}}$ \unboldmath \hspace{-3mm}
}}

\journalname{Numerische Mathematik}
\begin{document}

\title{Convergence of a finite volume scheme for the incompressible fluids
}


\author{S\'ebastien Zimmermann
}


\institute{S. Zimmermann \at
              17 rue Barr\`eme, 69006 Lyon - FRANCE\\
              Tel.: (+33)0472820337\\
              \email{Sebastien.Zimmermann@ec-lyon.fr}           
}

\date{Received: date / Revised: date}

\maketitle

\begin{abstract}
We consider a finite volume scheme for the two-dimensional
incompressible Navier-Stokes equations. We use a triangular mesh.
 The unknowns for the velocity and pressure are respectively piecewise constant and affine. We use a
  projection method to deal with the incompressibility constraint.
   The stability of the scheme has been proven in
\cite{zimm2}. We infer from it its  convergence.
 \subclass{Incompressible fluids \and Navier-Stokes equations \and
projection methods \and finite volume}
\end{abstract}

\section{Introduction}

We consider the flow of an incompressible fluid in a open bounded
polyhedral set $\Omega \subset \mathbb{R}^2$ during the time
interval $[0,T]$. The velocity field $\mathbf{u}:\Omega \times [0,T]
\to \mathbb{R}^2 $ and the pressure field $p:\Omega \times [0,T] \to
\mathbb{R}$ satisfy the Navier-Stokes equations
\begin{eqnarray}
&& \mathbf{u}_t - \frac{1}{\hbox{Re}} \, \mathbf{\Delta} \mathbf{u}
+(\U \cdot \Nabla) \U+\nabla p = \F \, ,
\label{eq:mom} \\
&& \hbox{div } \U=0 \, , \label{eq:incomp}
\end{eqnarray}
with the boundary and initial condition
\begin{equation*}
 \U|_{\partial \Omega}=0 \, , \hspace{2cm}  \U|_{t=0}=\U_0.
\end{equation*}
The terms $\mathbf{\Delta} \U$ and $(\U \cdot \Nabla)\U$ are
associated with the  physical phenomena of diffusion and convection,
respectively. The Reynolds number $\hbox{Re}$ measures the influence
of convection in the flow. For equations
(\ref{eq:mom})--(\ref{eq:incomp}), finite element and finite
difference methods are well known and mathematical studies are
available (see \cite{giraultr} for example). For finite volume
schemes, numerous computations have been conducted (\cite{kimchoi}
and \cite{boicaya} for example).
 However, few mathematical results are available in this case.  Let us cite {\sc Eymard and Herbin} \cite{herb3}
 and {\sc Eymard, Latch\'e and Herbin} \cite{eymard}.
 In order to deal with the incompressibility constraint, these works use a penalization  method.
 Another way is to use  the projection methods which have been
 introduced by {\sc Chorin} \cite{chorin} and {\sc Temam} \cite{temam}.
 This is the case in  {\sc Faure} \cite{faure} where
 the mesh is made of  squares.
 In {\sc Zimmermann} \cite{zimm1} the mesh is made of triangles, which allows more
 complex geometries. 
 In the present paper the mesh is also made of triangles, but we consider a different
 discretisation for the pressure. It leads to a linear system with a better-conditioned
 matrix.
The layout of the article is the following. We first introduce
(section \ref{subsec:defmaillage}) some notations and hypotheses on
the mesh. We define  (section \ref{subsec:espd}) the spaces we use
to approximate the velocity and pressure. We define also (section
\ref{subsec:opd}) the operators we use to approximate the
differential operators in (\ref{eq:mom})--(\ref{eq:incomp}).
 By combining this with a projection method, we build the scheme in section \ref{sec:presschema}.
  In order to provide a mathematical analysis, we
  state in section \ref{sec:propopd} that the differential operators
   in (\ref{eq:mom})--(\ref{eq:incomp}) and their discrete counterparts share similar properties.
 In particular, the discrete operators for the gradient and the
 divergence are adjoint. 
 We then prove in section \ref{sec:stab} the
convergence of the scheme.

\noindent We conclude with some notations. We denote by $\chi_I$ the characteristic function of
an interval  $I \subset \mathbb{R}$. We denote by ${\mathcal{C}}^\infty_0={\mathcal{C}}^\infty_0(\Omega)$
 the set of the functions with a compact support in $\Omega$.
The spaces $(L^2,|.|)$ and
$(L^\infty,\|.\|_\infty)$ are the usual Lebesgue spaces and we set
$L^2_0=\{q \in L^2 \, ; \int_\Omega q \, d\x=0 \}$. Their
vectorial counterparts are $(\L^2,|.|)$ and
$(\L^\infty,\|.\|_\infty)$ with $\L^2=(L^2)^2$ and
$\L^\infty=(L^\infty)$. For $k \in \mathbb{N}^*$,
$(H^k,\|\cdot\|_k)$ is the usual Sobolev space. Its vectorial
counterpart is $(\H^k,\|.\|_k)$ with $\H^k=(H^k)^2$. For $k=1$, the
functions of $\H^1$ with a null trace on the boundary form the space
$\H^1_0$. Also, we set $\Nabla \U=(\nabla u_1,\nabla u_2)^T$ if
$\U=(u_1,u_2) \in \H^1$.
If $\mathbf{X}\subset \L^2$ is a Banach space, we define
${\mathcal{C}}(0,T;\mathbf{X})$
 (resp. $L^2(0,T;\mathbf{X})$) as the set of the
 applications $\mathbf{g}:[0,T] \to \mathbf{X}$ such that
$t \to |\mathbf{g}(t)|$ is continous (resp. square integrable).
 The norm $\|.\|_{{\mathcal{C}}(0,T;\mathbf{X})}$
is defined by
$\|\mathbf{g}\|_{{\mathcal{C}}(0,T;\mathbf{X})}=\sup_{s\in[0,T]}
|\mathbf{g}(s)|$.
Finally in all calculations, $C$ is a generic positive constant, depending
only on $\Omega$, $\U_0$ and $\F$.
\section{Discrete setting}
\label{sec:defd}

First, we introduce the spaces and operators needed to build the mesh.

\subsection{The mesh}
\label{subsec:defmaillage}

Let ${\mathcal{T}}_h$ be a triangular mesh of $\Omega$:
$\overline{\Omega}=\cup_{K \in {\mathcal{T}}_h} K $. For each
triangle $K \in {\mathcal{T}}_h$, we denote by  $|K|$ its area and
${\mathcal{E}}_K$ the set of his edges.
 If $\sigma \in {\mathcal{E}}_K$, $\N_{K,\sigma}$  is the unit vector normal to $\sigma$ pointing outwards of $K$.

\noindent  The set of edges of the mesh is ${\mathcal{E}}_h=\cup_{K \in
{\mathcal{T}}_h} {\mathcal{E}}_K$. The length of an edge $\sigma \in
{\mathcal{E}}_h$ is $|\sigma|$. The set of edges inside $\Omega$
(resp. on the boundary) is ${\mathcal{E}}^{int}_h$
 (resp. ${\mathcal{E}}^{ext}_h$): ${\mathcal{E}}_h={\mathcal{E}}^{int}_h \cup {\mathcal{E}}^{ext}_h$.
 If $\sigma \in {\mathcal{E}}^{int}_h$, $K_\sigma$ and $L_\sigma$
 are the triangles sharing $\sigma$ as an edge.
If $\sigma \in {\mathcal{E}}^{ext}_h$, only the triangle $K_\sigma$
 inside $\Omega$ is defined.

\noindent  We denote by   $\x_K$  the circumcenter of a triangle $K$. We
assume that the measure of all interior angles of the triangles of
the mesh are below $\frac{\pi}{2}$, so that $\x_K \in K$. If $\sigma
\in {\mathcal{E}}^{int}_h$ (resp. $\sigma \in {\mathcal{E}}^{ext}_h$
) we set
 $d_\sigma=d(\x_{K_\sigma},\x_{L_\sigma})$ (resp.  $d_\sigma=d(\x_\sigma,\x_{K_\sigma})$).
We define for all edge $\sigma \in {\mathcal{E}}_h$:
  $\tau_\sigma=\frac{|\sigma|}{d_\sigma}$.
The maximum circumradius  of the triangles of the mesh is $h$. We
assume that there exists $C>0$  such that
\begin{equation*}
 \forall \, \sigma \in {\mathcal{E}}_h, \hspace{1cm} \; \; d(\x_{K_\sigma},\sigma) \ge C |\sigma|
\hspace{.5cm}   \hbox{and} \hspace{.5cm} |\sigma| \ge Ch.
\end{equation*}
It implies that there exists a constant $C>0$ such
that for all edge $\sigma \in {\mathcal{E}}_h$
\begin{equation}
\label{eq:proptaus}
  \tau_\sigma \ge C
\end{equation}
and for all triangles  $K \in {\mathcal{T}}_h$ we have (with $\sigma
\in {\mathcal{E}}_K$ and $h_{K,\sigma}$ the matching altitude)
\begin{equation}
\label{eq:propairet}
  |K|=\frac{1}{2} \, |\sigma| \, h_{K,\sigma} \ge \frac{1}{2} \, |\sigma| \, d(\x_K,\x_\sigma) \ge C \, h^2.
\end{equation}

\subsection{The discrete spaces}

\label{subsec:espd}

\noindent   We first define
 \begin{equation*}
 \label{eq:defp0}
   P_0= \{ q \in L^2 \; ; \; \forall \, K \in {\mathcal{T}}_h, \; \; q|_K \hbox{ is a constant} \} \, ,
   \hspace{1cm} \P_0=(P_0)^2.
\end{equation*}
For the sake of concision, we set for all $q_h \in P_0$ (resp. $\V_h
\in \P_0$) and all triangle $K \in {\mathcal{T}}_h$:
 $q_K=q_h|_K$ (resp. $\V_K=\V_h|_K$).
Although $\P_0 \not \subset \H^1$, we define the discrete equivalent
of a $\H^1$ norm as follow. For all $\V_h \in \P_0$ we set
\begin{equation}
\label{eq:defh1d}
  \|\V_h\|_h=\left( \sum_{\sigma \in {\mathcal{E}}^{int}_h} \tau_{\sigma} \, |\V_{L_\sigma} - \V_{K_\sigma}|^2
  +\sum_{\sigma \in {\mathcal{E}}^{ext}_h} \tau_\sigma \, |\V_{K_\sigma}|^2 \right)^{1/2}.
\end{equation}
 We have \cite{eymgal} a discrete Poincar\'e inequality for $\P_0$: there exists $C>0$ such that for all
  $\V_h \in \P_0$
 \begin{equation}
 \label{eq:inpoinp0}
  |\V_h| \le C \, \|\V_h\|_h.
 \end{equation}
From the norm $\|.\|_h$ we deduce a dual norm.
 For all $\V_h \in \P_0$ we set
\begin{equation}
\label{eq:defnormdp0}
  \|\V_h\|_{-1,h}=\sup_{\psig_h \in \P_0} \frac{(\V_h,\psig_h)}{\|\psig_h\|_h}.
\end{equation}
For all $\U_h \in \P_0$ and $\V_h \in \P_0$ we have
  $(\U_h,\V_h) \le \|\U_h\|_{-1,h} \, \|\V_h\|_h$.
 We define the projection operator  $\Pi_{\P_0}: \L^2 \to \P_0$ as follows.
 For all $\W \in \L^2$, $\Pi_{\P_0} \W \in \P_0$ is given by
\begin{equation}
\label{eq:defprojp0}
   \forall \, K \in {\mathcal{T}}_h \, ,  \hspace{1cm} (\Pi_{\P_0} \W)|_K=\frac{1}{|K|} \int_K \W(\x) \, d\x.
\end{equation}
We easily check that for all  $\W \in \L^2$ and $\V_h \in \P_0$ we
have $(\Pi_{\P_0} \W, \V_h)=(\W,\V_h)$. We deduce from this that
$\Pi_{\P_0}$ is stable for the $\L^2$ norm.
 We define also the operator $\widetilde \Pi_{\P_0}:\H^2 \to \P_0$.
 For all $\V\in \H^2$, $\widetilde \Pi_{\P_0} \V \in \P_0$
 is given by
\begin{equation*}
  \forall \, K \in {\mathcal{T}}_h \, , \hspace{1cm} \widetilde \Pi_{\P_0} \V|_K=\V(\x_K).
\end{equation*}
According to the Sobolev embedding theorem, $\V \in \H^2$ is a.e. equal
to a continuous function. Therefore the definition above makes
sense. 
One checks \cite{zimm1} that there exists $C>0$ such that
\begin{equation}
\label{eq:esterrintp0}
|\V-\widetilde \Pi_{\P_0}\V| \le C \, h\, \|\V\|_2
\end{equation}
for all $\V \in \H^1$.
\noindent  We introduce also the finite element spaces
\begin{eqnarray}
\label{eq:defp1nc}
  P^d_1 &=& \{ v \in L^2 \; ; \; \forall \, K \in {\mathcal{T}}_h, \; \; v|_K \hbox{ is affine} \} \, , \nonumber\\
  P^{nc}_1 &=& \{ v_h \in P^d_1 \; ; \; \forall \, \sigma \in {\mathcal{E}}^{int}_h, \,
  v_h|_{K_\sigma}(\x_\sigma)=v_h|_{L_\sigma}(\x_\sigma) \, , \quad \P^{nc}_1=(P^{nc}_1)^2. \nonumber 
 \end{eqnarray}
If $q_h \in P^{nc}_1$, we
have usually $\nabla q_h \not \in \L^2$. Therefore we define the operator
$\nabla_h: P^{nc}_1 \hspace{-.2cm}\to \P_0 $
 by setting  for all $q_h \in P_0$ and all triangle  $K \in {\mathcal{T}}_h$
\begin{equation}
\label{eq:defg}
  \nabla_h q_h |_K = \frac{1}{|K|} \, \int_K \nabla q_h \, d\x.
\end{equation}
We define the projection operator $\Pi_{P^{nc}_1}$. For all $q \in
H^1$, $\Pi_{P^{nc}_1} q$ is given by
\begin{equation}
\label{eq:defpp1nc}
   \forall \, \sigma \in {\mathcal{E}}_h \, , \hspace{1cm}
    \int_\sigma (\Pi_{P^{nc}_1} q) \, d\sigma=\int_\sigma q \, d\sigma.
\end{equation}
We also set $\Pi_{\P^{nc}_1}=(\Pi_{P^{nc}_1})^2$. One checks that there exists $C>0$ such that
\begin{equation}
\label{eq:esterrintp1nc}
  \left|\nabla q-\nabla_h(\Pi_{P^{nc}_1}q)\right| \le C \, h\, \|q\|_2 \, , \hspace{1cm} |\V-\Pi_{\P^{nc}_1} \V| \le C \, h\, \|\V\|_1 \, ,
\end{equation}
for all $q \in H^1$ and $\V \in \H^1$.

\noindent  We also use the Raviart-Thomas spaces
\begin{eqnarray*}
\label{eq:defrt0}
  \mathbf{RT^d_0}&=& \{ \V_h \in \P^d_1 \; ; \; \forall \, \sigma \in {\mathcal{E}}_K,
   \; \V_h|_K \cdot \N_{K,\sigma} \hbox{ is constant,}  \; \hbox{ and }  \;
    \V_h \cdot \N|_{\partial \Omega}=0  \} \, , \nonumber \\
  \mathbf{RT_0}&=& \{ \V_h \in \mathbf{RT^d_0} \; ; \; \forall \, K \in {\mathcal{T}}_h,
   \; \forall \, \sigma \in {\mathcal{E}}_K,
 \; \;    \V_h|_{K_\sigma} \cdot \N_{K_\sigma,\sigma} = \V_h|_{L_\sigma} \cdot \N_{K_\sigma,\sigma}\}.
\end{eqnarray*}
For all $\V_h \in \mathbf{RT_0}$, $K \in {\mathcal{T}}_h$ and
$\sigma \in {\mathcal{E}}_K$ we set
  $(\V_h \cdot \N_{K,\sigma})_\sigma=\V_h|_K \cdot \N_{K,\sigma}$.
\noindent We define the operator  $\Pi_{\mathbf{RT_0}}:\H^1 \to
\mathbf{RT_0}$. For all $\V \in \H^1$, $\Pi_{\mathbf{RT_0}} \V \in
\mathbf{RT_0}$ is given by
\begin{equation}
\label{eq:defprt0} \forall \, K \in {\mathcal{T}}_h \, ,
\hspace{.5cm}
  \forall \, \sigma\in {\mathcal{E}}_K \, , \hspace{1cm}
    (\Pi_{\mathbf{RT_0}} \V \cdot \N_{K,\sigma})_\sigma=\frac{1}{|\sigma|}\int_\sigma \V \, d\sigma.
\end{equation}
One checks \cite{brezzfor} that there exists a
constant $C>0$ such that for all $\V \in \H^1$
\begin{equation}
\label{eq:esterrirt0}
  |\V-\Pi_{\mathbf{RT_0}} \V| \le C \, h \, \|\V\|_1.
\end{equation}

\subsection{The discrete operators}

\label{subsec:opd}

The equations (\ref{eq:mom})--(\ref{eq:incomp}) use the
differential operators gradient, divergence and laplacian.
 Using the spaces of section \ref{subsec:espd}, we now define their discrete counterparts.
The discrete gradient $ \nabla_h: P^{nc}_1 \to \P_0$ is
defined by (\ref{eq:defg}). \noindent  The discrete divergence
operator $\hbox{div}_h: \P_0 \to P^{nc}_1$ is built so that it is
adjoint to the operator $\nabla_h$ (proposition \ref{prop:propadjh} below). We set for all $\V_h \in \P_0$
and all triangle $K \in {\mathcal{T}}_h$
 \begin{eqnarray}
\label{eq:defdivh1} &&  \forall \, \sigma \in {\mathcal{E}}^{int}_h,
\hspace{2cm}  (\hbox{div}_h \, \V_h) (\x_\sigma)=
\frac{3 \, |\sigma|}{|K_\sigma|+|L_\sigma|} \, (\V_{L_\sigma} - \V_{K_\sigma}) \cdot \N_{K,\sigma} \; ; \nonumber \\
&&  \forall \, \sigma \in {\mathcal{E}}^{ext}_h, \hspace{2cm}
(\hbox{div}_h \, \V_h) (\x_\sigma)= -\frac{3 \,
|\sigma|}{|K_\sigma|+|L_\sigma|} \, \V_{K_\sigma} \cdot
\N_{K,\sigma}. \label{eq:defdivh2}
\end{eqnarray}
The first discrete laplacian $\Delta_h:P^{nc}_1 \to P^{nc}_1$ is given by
\begin{equation*}
\label{eq:deflap} \forall \, q_h \in P^{nc}_1 \, , \hspace{1cm} \Delta_h q_h=\hbox{div}_h (\nabla_h q_h).
\end{equation*}
\noindent The second discrete laplacian
$ \Deltat_h: \P_0 \to \P_0$ is the usual operator in finite volume
schemes \cite{eymgal}. We set for all $\V_h \in \P_0$ and all triangle $K
\in {\mathcal{T}}_h$
\begin{equation}
\label{eq:defl}
  \Deltat_h \V_h|_K
  = \frac{1}{|K|} \sum_{\sigma \in {\mathcal{E}}_K \cap {\mathcal{E}}^{int}_h}
 \tau_\sigma \, (\V_{L_\sigma} - \V_{K_\sigma} )
-\frac{1}{|K|}\sum_{\sigma \in {\mathcal{E}}_K \cap
{\mathcal{E}}^{ext}_h}
 \tau_\sigma \, \V_{K_\sigma}.
\end{equation}

\noindent In order to approximate the convection term  $(\U \cdot \Nabla) \U$
in (\ref{eq:mom}), we define a bilinear form $\bt_h: \mathbf{RT_0}
\times \P_0 \to \P_0$ using the well-known \cite{eymgal} upwind scheme. For all
$\U_h \in \P_0$, $\V_h \in \P_0$, and all triangle $K \in
{\mathcal{T}}_h$ we set
\begin{equation*}
\label{eq:defbth}
  \bt_h(\U_h,\V_h)\big|_K=
\frac{1}{|K|} \sum_{\sigma \in {\mathcal{E}}_K \cap
{\mathcal{E}}^{int}_h} \,
   |\sigma| \, \Big( (\U \cdot \N_{K,\sigma})^+_\sigma \, \V_K +
(\U \cdot \N_{K,\sigma})^-_\sigma \, \V_{L_\sigma}  \Big).
\end{equation*}
We have set $a^+=\max(a,0)$, $a^-=\min(a,0)$ for all
$a\in\mathbb{R}$.
Lastly, we define the trilinear form $ \b_h: \mathbf{RT_0}
\times \P_0 \times \P_0 \to \mathbb{R}$ as follows. For all $\U_h
\in \mathbf{RT_0}$, $\V_h \in \P_0$, $\W_h \in \P_0$, we set
\begin{equation*}
\label{eq:defbh}
  \b_h(\U_h,\V_h,\W_h)=\sum_{K \in {\mathcal{T}}_h} |K| \, \W_K \cdot
   \bt_h(\U_h,\V_h)\big|_K .
\end{equation*}

\section{The scheme}
\label{sec:presschema}

We have defined in section \ref{sec:defd} the discretization in space. We now have to define the discretization in time,
 and treat  the incompressibility constraint (\ref{eq:incomp}).
We use a projection method to this end. This kind
of method has been introduced by {\sc Chorin} \cite{chorin}
 and {\sc Temam} \cite{temam}.
The time interval $[0,T]$ is split with a time step $k$:
$[0,T]=\bigcup_{n=0}^N [t_n,t_{n+1}]$ with $N \in \mathbb{N}^*$ et
$t_n=n \, k$ for all $n\in\{0,\dots,N\}$.
We start with the initial values
\begin{equation*}
  \U^0_h \in \P_0 \cap \mathbf{RT_0} \, , \hspace{1cm}   \U^1_h \in \P_0 \cap \mathbf{RT_0} \, ,
\hspace{1cm}   p^1_h \in P^{nc}_1 \cap L^2_0.
\end{equation*}
For all  $n \in \{1,\dots,N\}$, $(\Ut^{n+1}_h,p^{n+1}_h,\U^{n+1}_h)$
is deduced from $(\Ut^n_h,p^n_h,\U^n_h)$ as follows. 

\begin{itemize}
\item $\Ut^{n+1}_h \in \P_0$ is given by
\begin{align}
\label{eq:mombdf} \frac{3 \, \Ut^{n+1}_h -4 \, \U^n_h+\U^{n-1}_h}{2
\, k} &- \frac{1}{\hbox{Re}} \, \Deltat_h \Ut^{n+1}_h \nonumber \\
&+\bt_h(2\, \U^n_h-\U^{n-1}_h,
\Ut^{n+1}_h)+\nabla_h p^n_h =\F^{n+1}_h  \, ,
\end{align}

\item $p^{n+1}_h \in P^{nc}_1 \cap L^2_0$ is the solution of
\begin{equation*}
\label{eq:projia} \Delta_h (p^{n+1}_h-p^n_h)=\frac{3}{2 \, k}\,
\hbox{div}_h \, \Ut^{n+1}_h  \, ,
\end{equation*}

\item $\U^{n+1}_h \in  \P_0$ is given by
\begin{equation}
\label{eq:projib}
 \U^{n+1}_h = \Ut^{n+1}_h -\frac{2\, k}{3} \, \nabla_h (p^{n+1}_h - p^n_h).
\end{equation}
\end{itemize}
We have proven in \cite{zimm1} that the scheme is well defined. In particular the term
$\bt_h(2\, \U^n_h-\U^{n-1}_h)$ in (\ref{eq:mombdf}) is defined thanks to the following result.
\begin{proposition}
\label{prop:umrt0} For $m\in\{0,\dots,N\}$ we have  $\U^m_h \in
\mathbf{RT_0}$.
\end{proposition}
Note also that  for $m\in\{0,\dots,N\}$ we have $\hbox{div} \, \U^m_h=0$, since $\U^m_h \in \P_0$. Thus the incompressibility
condition (\ref{eq:incomp}) is fullfilled.



\section{Properties of the discrete operators}
\label{sec:propopd}

The operators defined in section \ref{subsec:opd} have the following
properties \cite{zimm1}.
\begin{proposition}
\label{prop:stabbth} There exists a constant $C>0$ such that for all
$\U_h \in \mathbf{RT_0}$ satisfying   $\hbox{ \rm div} \, \U_h=0$,
$\V_h \in \P_0$, $\W_h\in \P_0$:
\begin{equation*}
\label{eq:majbh} |\b_h(\U_h,\V_h,\V_h)| \le C \, |\U_h| \,
\|\V_h\|_{h} \, \|\W_h\|_h.
\end{equation*}
\end{proposition}
\begin{proposition}
\label{prop:propadjh} For all $\V_h \in \P_0$ and $q_h \in
P^{nc}_1$: $(\V_h, \nabla_h q_h)=-(q_h,\hbox{\rm div}_h \, \V_h)$.
\end{proposition}

\begin{proposition}
\label{prop:coerlap} For all $\U_h \in \P_0$ and $\V_h \in \P_0$:
$-(\Deltat_h \U_h,\U_h)=\|\U_h\|^2_h \label{eq:propdh2}$ and
$-(\Deltat_h \U_h,\V_h) \le \|\U_h\|_h \,
\|\V_h\|_h\label{eq:propdh1}$.

\end{proposition}
If $\V \in \H^1$ we have $|\hbox{div} \, \V| \le \|\V\|_1$. The
operator $\hbox{div}_h$ has a similar property.
\begin{proposition}
  \label{propstabdivh}
  There exists a constant  $C>0$ such that for all $\V_h \in \P_0$
  \begin{equation*}
  \label{eq:stabdivh}
    |\hbox{\rm div}_h \, \V_h| \le C \, \|\V_h\|_h.
  \end{equation*}
\end{proposition}
\noindent {\sc Proof.} Using a quadrature formula we have
\begin{equation*}
  |\hbox{div}_h \, \V_h|^2=\sum_{K \in {\cal{T}}_h} \frac{|K|}{3}  \sum_{\sigma \in {\cal{E}}_K}
  \left|(\hbox{div}_h \, \V_h)(\x_\sigma)\right|^2.
\end{equation*}
Let $K \in {\cal{T}}_h$.
Using definition and   (\ref{eq:propairet}) we have
\begin{eqnarray*}
&&\forall \, \sigma \in {\cal{E}}_K \cap {\cal{E}}^{int}_h \, ,
\hspace{1cm}  \left|(\hbox{div}_h \, \V_h)(\x_\sigma)\right|^2  \le
\frac{C}{|K|} \, |\V_{L_\sigma}-\V_K|^2 \, ; \\
&&\forall \, \sigma \in {\cal{E}}_K \cap {\cal{E}}^{ext}_h \, ,
\hspace{1cm}  \left|(\hbox{div}_h \, \V_h)(\x_\sigma)\right|^2  \le
\frac{C}{|K|} \, |\V_K|^2.
\end{eqnarray*}
Thus: 
   $|\hbox{div}_h \, \V_h|^2 \le C \sum_{K \in {\cal{T}}_h} \! \! \! \left( \sum_{\sigma \in {\cal{E}}_K
    \cap {\cal{E}}^{int}_h} |\V_{L_\sigma}-\V_K|^2
+\sum_{\sigma \in {\cal{E}}_K   \cap {\cal{E}}^{ext}_h} \! |\V_K|^2
\right)$.
Writing the sum over the triangles  as a sum over the edges, we get
\begin{equation*}
\hspace{.45cm}  |\hbox{div}_h \, \V_h|^2 \le C \, \left(\sum_{\sigma
\in {\cal{E}}^{int}_h} \tau_\sigma \, |\V_{L_\sigma}-\V_K|^2
 +\sum_{\sigma \in {\cal{E}}^{ext}_h} \tau_\sigma \, |\V_K|^2 \right)\le C \, \|\V_h\|^2_h. \hspace{.45cm} \qed
\end{equation*}

\begin{proposition}
\label{prop:adjlap} If $\U_h \in \P_0$ and $\V_h \in \P_0$ we have
  $(\Deltat_h \U_h, \V_h)=(\U_h,\Deltat_h \V_h)$.
\end{proposition}
{\sc Proof.} Using definition (\ref{eq:deflap}) one checks that
\begin{eqnarray*}
  (\Deltat_h \U_h, \V_h)
&=& \sum_{\sigma \in {\cal{E}}^{int}_h}
 \tau_\sigma \, (\V_{L_\sigma}-\V_{K_\sigma}) \cdot (\U_{L_\sigma} -\U_{K_\sigma})
-\sum_{\sigma \in {\cal{E}}^{ext}_h}  \tau_\sigma \, \V_{K_\sigma}
\cdot \U_{K_\sigma} \\
&=&(\Deltat_h \V_h, \U_h). \hspace{7.3cm} \qed
\end{eqnarray*}

\begin{proposition}
\label{propconslap} There exists $C>0$ such that for all $\V \in
\H^2$ satisfying  $\Nabla \V \cdot \N|_{\partial \Omega}=\mathbf{0}$
\begin{equation*}
 \label{eq:conslap} \|\Pi_{\P_0} (\mathbf{\Delta} \V) -
\mathbf{\Delta}_h(\widetilde \Pi_{\P_0} \V)\|_{-1,h} \le C \, h \,
\|\V\|_2.
\end{equation*}
\end{proposition}
\noindent {\sc Proof.} Let $\psig_h \in \P_0$. We have
\begin{equation*}
\big( \Pi_{\P_0} (\mathbf{\Delta} q) - \mathbf{\Delta}_h(\widetilde
\Pi_{\P_0} \V),\psig_h \big) =\sum_{K \in {\cal{T}}_h} |K| \, \left.
\big( \Pi_{\P_0} (\mathbf{\Delta} q) - \mathbf{\Delta}_h(\widetilde
\Pi_{\P_0} \V)\big)\right|_K \cdot \psig_K.
\end{equation*}
 For all $K \in {\cal{T}}_h$, using (\ref{eq:deflap}) and
the divergence formula,  we get
\begin{eqnarray*}
&&  |K|\, \left. \big( \Pi_{\P_0} (\mathbf{\Delta} \V) - \mathbf{\Delta}_h(\widetilde \Pi_{\P_0} \V) \big) \, \right|_K \nonumber \\
&&= \sum_{\sigma \in {\cal{E}}_K \cap {\cal{E}}^{int}_h} \left(
\int_\sigma \Nabla \V \cdot \N_{K,\sigma} \, d\sigma -
\tau_{\sigma} \, \big(\V(\x_{L_\sigma})-\V(\x_{K})\big)\right).
\end{eqnarray*}
Thus, by writing the sum over the triangles as a sum over the edges, we get
\begin{equation*}
\big( \Pi_{\P_0} (\mathbf{\Delta} \V) - \mathbf{\Delta}_h(\widetilde \Pi_{\P_0} \V),\psig_h \big)
=\sum_{\sigma \in {\cal{E}}^{int}_h}
(\psig_{L_\sigma}-\psig_{K_\sigma}) \, \R_\sigma
\end{equation*}
with
  $\R_\sigma=
\int_\sigma \Big( \Nabla \V \cdot \N_{K_\sigma,\sigma}  -
\frac{1}{d_\sigma}  \, \big( \V(\x_{L_\sigma})-\V(\x_{K_\sigma})
\big) \Big) \, d\sigma$.
We denote by $D_\sigma$ the quadrilatere defined by $\x_{K_\sigma}$, $\x_{L_\sigma}$ and the endpoints of $\sigma$. Using a Taylor
 expansion and a density argument, we get as in \cite{eymgal} that
\begin{equation*}
\label{eq:conslap7}
   |R_\sigma|
\le  C \, h  \left( \int_{D_\sigma} \sum_{i=1}^2|\H(v_i)(\mathbf{y})|^2 \,
d\mathbf{y} \right)^{1/2}.
\end{equation*}
Thus, using the Cauchy-Schwarz inequality, we get
\begin{eqnarray*}
&& |\big( \Pi_{\P_0} (\mathbf{\Delta} q) - \mathbf{\Delta}_h(\widetilde \Pi_{\P_0} \V),\psig_h \big)| \\
&& \le C \, h \left(\sum_{\sigma \in {\cal{E}}^{int}_h}
|\psig_{L_\sigma}-\psig_{K_\sigma}|^2 \right)^{1/2}
\left(\sum_{i=1}^2 \sum_{\sigma \in {\cal{E}}^{int}_h}
\int_{D_\sigma}
 |\H(v_i)(\mathbf{y})|^2
 \, d\mathbf{y} \right)^{1/2}.
\end{eqnarray*}
According to  (\ref{eq:proptaus})
\begin{equation*}
\sum_{\sigma \in {\cal{E}}^{int}_h}
|\psig_{L_\sigma}-\psig_{K_\sigma}|^2  \le C
\sum_{\sigma \in {\cal{E}}^{int}_h} \tau_\sigma \,
|\psig_{L_\sigma}-\psig_{K_\sigma}|^2  \le C \, \|\psig\|^2_h.
\end{equation*}
Thus
$\left|\big( \Pi_{\P_0} (\mathbf{\Delta} \V) -
\mathbf{\Delta}_h(\widetilde \Pi_{\P_0} \V),\psig_h \big)\right| \le
C \, h \,  \|\psig_h\|_h \, \|\V\|_2$.
Using (\ref{eq:defnormdp0}) we get the result. \qed

\section{Convergence of the scheme}
\label{sec:stab}

We first recall the stability result that has been proven in \cite{zimm2}.
We deduce from it  an estimate on the Fourier transform of the
computed velocity (lemma \ref{lemcompvi}). Using a result on space
$\P_0$, we infer from it the convergence of the scheme (theorem
\ref{theo:conv}).
\noindent One shows that if the data  $\U_0$ et $\F$ fulfill
 a compatibility condition   \cite{heywood},
there exists a  solution  $(\overline \U,\overline p)$ to equations
(\ref{eq:mom})--(\ref{eq:incomp}) such that
 $\overline \U \in {\cal{C}}(0,T;\H^2)$   ,
 $ \nabla \overline p\in {\cal{C}}(0,T;\L^2)$.
We assume from now on that there exists $C>0$ such that
\begin{equation*}
  {\bf (HI)}    \hspace{.7cm}  |\U^0_h-\U_0| + \frac{1}{h} \, \|\U^1_h-\U(t_1)\|_{\infty}
      +|p^1_h-p(t_1)| \le C \, h \, , \hspace{.7cm} |\U^1_h-\U^0_h|\le C \, k.
\end{equation*}
Let us recall the following result \cite{zimm2}.
\begin{theorem}
\label{theo:stab} We assume that the initial values of the scheme
fulfill  {\bf (HI)}. There exists a constant  $C>0$ such that for
all $m \in \{2,\dots,N\}$
\begin{equation*}
\label{eq:eststabv}
  |\U^m_h|+k\sum_{n=2}^m \|\Ut^n_h\|^2_h+\frac{1}{k} \, |\U^m_h-\U^{m-1}_h|
+\frac{1}{k}\sum_{n=2}^m \|\Ut^n_h-\Ut^{n-1}_h\|^2_h
   \le C
\end{equation*}
and
\begin{equation*}
\label{eq:eststabp}
  k \sum_{n=2}^m |p^n_h|^2 +k \, |\nabla_h p^m_h|+|\nabla_h(p^m_h-p^{m-1}_h)|\le C.
\end{equation*}
\end{theorem}
From now on we set $\Ut^1_h=\U^1_h$ for the sake of conveniance. One deduces from hypothesis {\bf (HI)} \cite{zimm1}
that $|p^1_h| \le C$ and $\|\Ut^1_h\|_h \le C.$
Now, let $\eps=\max(h,k)$. We study the behaviour of the scheme as $\eps \to 0$. We
define the applications
$\U_{\eps}:\mathbb{R} \to \P_0$ ,  $\Ut_{\eps}:\mathbb{R} \to \P_0$
, $\Ut^{c}_{\eps}:\mathbb{R} \to \P_0$, $p_{\eps}:\mathbb{R} \to
P^{nc}_1$ and  $\F_{\eps}:\mathbb{R} \to \P_0$
as follows. For all $n \in \{0,\dots,N-1\}$ and all  $t
\in[t_n,t_{n+1}]$ we set
\begin{eqnarray*}
\label{eq:defueps} &&   \U_{\eps}(t)  = \U^n_h \; , \hspace{.5cm}
 \Ut_{\eps}(t)  = \Ut^{n+1}_h, \hspace{.5cm}
\Ut^c_{\eps}(t)=\Ut^n_h + \frac{1}{k} \, (t-t_n) \, (\Ut^{n+1}_h - \Ut^n_h) \, , \\
&&   p_\eps(t)=p^n_h, \hspace{.5cm} \F_\eps(t)=\Pi_{\P_0}
\F(t_{n+1}) \, , \nonumber
\end{eqnarray*}
\noindent and for all $t \not \in [0,T]$ we set
   $\U_{\eps}(t)  = \Ut_{\eps}(t)=\Ut^c_{\eps}(t)=\F_\eps(t)=\mathbf{0}$, $p_\eps(t)=0$.
We recall that the Fourier transform $\widehat \V$ of  a function
$\V \in \L^1(\mathbb{R})$  is defined by
\begin{equation}
\label{eq:deftf} \forall \, \tau \in \mathbb{R}, \; \; \; \; \;
\widehat \V (\tau) = \int_{\mathbb{R}} e^{-2i\pi\tau t} \,  \V(t) \,
dt.
\end{equation}
We have the following result.
\begin{lemma}
\label{lemcompvi}

Let $0<\gamma <\frac{1}{4}$. There exists  $C>0$ such that for all
$\eps>0$
\begin{equation*}
\label{eq:estcomp}
  \int_{\mathbb{R}} |\tau|^{2\gamma} \, |\widehat {\Ut}_{\eps}(\tau)|^2 \, d\tau \le C.
\end{equation*}
\end{lemma}
\noindent {\sc Proof.} Let $\chi_I$ be the characteristic function
of an interval  $I \subset \mathbb{R}$. We define the application
$\G_{\eps}: \mathbb{R} \to \P_0$ as follows. For all $t\not \in
[t_1,T]$ we set $\G_\eps(t)=\mathbf{0}$.  For all $t \in [t_1,T]$,
 $\G_\eps(t) \in \P_0$ is the solution of
\begin{equation*}
\label{eq:convgi}
 \Deltat_h \G_{\eps}  =
 \Deltat_h \Ut_{\eps}+\F_{\eps}-\bt_h\big(2 \, \U_\eps-\U_\eps(t-k),\Ut_\eps\big)
 +\P_\eps
\end{equation*}
with
$\P_\eps= -  \nabla_h p_{\eps}-2 \, \frac{\Ut_{\eps}(t-k)-
\U_{\eps}}{k} +\frac{\Ut_{\eps}(t-2 \, k)- \U_{\eps}(t-k)}{2 \, k}
\, \chi_{[t_2,T]}$.
We have omitted most of the time dependancies for the sake of
concision. Let us estimate $\G_{\eps}$. We have
\begin{align}
 -(\Deltat_h \G_{\eps}, \G_{\eps})
&= -(\Deltat_h \Ut_{\eps}, \G_{\eps})
 -(\F_{\eps},\G_{\eps}) \label{eq:conveqgeps} \\
 &+\b_h\big(2 \, \U_\eps-\U_\eps(t-k),\Ut_\eps,\G_\eps\big)
+(\P_\eps,\G_\eps). \nonumber
\end{align}
 According to proposition \ref{prop:coerlap} we have
\begin{equation*}
\label{eq:conveq7} -(\Deltat_h \G_{\eps},
\G_{\eps})=\|\G_{\eps}\|^2_h \, ,
\hspace{1cm}  -(\Deltat_h \Ut_{\eps}, \G_{\eps}) \le
\|\Ut_{\eps}\|_h \,  \|\G_{\eps}\|_h.
\end{equation*}
Using the Cauchy-Schwarz inequality and  (\ref{eq:inpoinp0}) we have
\begin{equation*}
-(\F_{\eps}, \G_{\eps}) \le |\F_{\eps}| \, |\G_{\eps}| \le C \,
|\F_{\eps}| \, \|\G_{\eps}\|_h. \label{eq:conveq81}
\end{equation*}
According to  proposition \ref{prop:stabbth} and theorem
\ref{theo:stab}
\begin{eqnarray*}
\label{eq:conveq82} && \left|\b_h\big(2 \,
\U_\eps-\U_\eps(t-k),\Ut_\eps,\G_\eps\big)\right| \\
&&\le C \, |2 \, \U_\eps-\U_\eps(t-k)| \, \|\Ut_\eps\|_h \,
\|\G_\eps\|_h
 \le  C  \, \|\Ut_\eps\|_h \, \|\G_\eps\|_h.
\end{eqnarray*}
Using (\ref{eq:projib}) we have
\begin{equation*}
\;  \P_\eps =-\left(1+\frac{4}{3} \, \chi_{[t_3,T]}\right)
\nabla_h p_\eps  +\frac{5}{3} \, \chi_{[t_3,T]} \, \nabla_h
p_\eps(t-k) -\frac{1}{3} \, \chi_{[t_3,T]} \, \nabla_h p_\eps(t-2 \,
k).
\end{equation*}
Using proposition \ref{prop:propadjh} and the Cauchy-Schwarz
inequality we get
\begin{equation*}
 \left|(\P_\eps,\G_\eps)\right| \le C \, \Big( |p_\eps|+ \chi_{[t_3,T]} \, |p_\eps(t-k)|+
\chi_{[t_3,T]} \, |p_\eps(t-2 \, k)| \Big) \, |\hbox{div}_h \,
\G_\eps|.
\end{equation*}
Using proposition \ref{propstabdivh} we have
\begin{equation*}
 \left|(\P_\eps,\G_\eps)\right| \le C \, \Big( |p_\eps|+ \chi_{[t_3,T]} \, |p_\eps(t-k)|+
\chi_{[t_3,T]} \, |p_\eps(t-2 \, k)| \Big) \, \|\G_\eps\|_h.
\end{equation*}
Let us plug these estimates into (\ref{eq:conveqgeps}). By simplifying
by $\|\G_{\eps}\|_h$ and integrating from $t=t_1$ to $T$ we get
\begin{equation*}
\label{eq:conveq13} \int^T_{t_1} \|\G_{\eps}\|_h \, dt \le C+C
\left( \int^T_{t_1}   |p_{\eps}| \, dt+\int^T_{t_1} |\F_{\eps}| \,
dt + \int^T_{t_1}  \|\Ut_{\eps}\|_h \, dt \right).
\end{equation*}
According to the  Cauchy-Schwarz inequality and theorem \ref{theo:stab}
 \begin{equation*}
 \label{eq:conveq14}
  \int^T_0 |p_{\eps}(t)| \, dt \le
\sqrt{T}
 \left(  \int^T_{t_1} |p_{\eps}(t)|^2 \, dt \right)^{1/2}
 \le \sqrt{T}   \left(  k\sum_{n=1}^N |p^n_h|^2 \right)^{1/2}
\le C.
 \end{equation*}
Thanks to the stability of  $\Pi_{\P_0}$ for the  $\L^2$ norm we
have
\begin{equation*}
\label{eq:conveq15} \int^T_0 |\F_{\eps}(t)| \, dt =k\sum^N_{n=1}
|\Pi_{\P_0} \F(t_{n})|
 \le k\sum^N_{n=1} |\F(t_{n})|
 \le k\sum^N_{n=1} \|\F\|_{{\cal{C}}(0,T;\L^2)} \le C.
\end{equation*}
And thanks to the Cauchy-Schwarz inequality and theorem
\ref{theo:stab}
\begin{equation*}
\label{eq:conveq6} \int^T_{t_1} \|\Ut_{\eps}(t)\|_h \, dt \le
\sqrt{T} \int^T_{t_1} \|\Ut_{\eps}(t)\|^2_h \, dt \le C \, k
\sum^N_{n=2} \|\Ut^n_h\|^2_h \le C.
\end{equation*}
Thus, since $\G_\eps(t)=0$ for $t \in [0,t_1]$, we get
 $\int^T_0 \|\G_{\eps}(t)\|_h=\int^T_{t_1} \|\G_{\eps}(t)\|_h  \, dt \le C$.
Using definition (\ref{eq:deftf}) we obtain finally
\begin{equation}
\label{eq:conveq135}
 \forall \, \tau \in \mathbb{R} \,, \hspace{1cm} \|\widehat \G_{\eps}(\tau)\|_h \le C.
\end{equation}
With this estimate we can now prove the result. Since the function
$\Ut^c_{\eps}$ is piecewise ${\cal{C}}^1$ on $\mathbb{R}$, and
discontinous for $t=0$ and $t=T$, equation (\ref{eq:mombdf}) reads
\begin{eqnarray*}
\label{eq:eqconv14} \frac{3}{2} \, \frac{d \tilde \U^c_{\eps}}{dt}
-\frac{1}{2} \, \frac{d \tilde \U^c_{\eps}}{dt}(t-k) &=&
\mathbf{\Delta}_h \G_{\eps}+\frac{3}{2} \, (\Ut^0_h \, \delta_0
-\Ut^N_h \, \delta_T)-\frac{1}{2} \,
  (\Ut^1_h \, \delta_{t_1}-\Ut^N_h \, \delta_{T+k}) \\
  &+& \frac{3}{2} \, \frac{\Ut^1_h-\Ut^0_h}{k} \, \chi_{[0,t_1]}-\frac{1}{2} \, \frac{\Ut^N_h-\Ut^{N-1}_h}{k} \, \chi_{[T,T+k]}
\end{eqnarray*}
where $\delta_0$, $\delta_{t_1}$, $\delta_T$ and $\delta_{T+k}$ are
Dirac distributions located respectively in $0$, $t_1$, $T$ and
$T+k$.
 Let $\tau \in \mathbb{R}$. Applying the Fourier transform we get
\begin{equation*}
\label{eq:eqconv15}
   -2  i  \pi  \tau \, \left( \frac{3}{2} -\frac{1}{2} \, e^{-2 i\pi \tau k}\right) \widehat \Ut_{\eps}(\tau)  \\
 =    \mathbf{\Delta}_h \widehat \G_{\eps}(\tau) +\alphag
 \end{equation*}
 with
 \begin{eqnarray*}
\alphag&=&  \left( \frac{3}{2} \, (\Ut^0_h-\Ut^N_h \, e^{-2 i \pi
T})
-\frac{1}{2} \, (\Ut^1_h-\Ut^N_h \, e^{-2 i \pi T})   \right) \frac{e^{-2i\pi k}-1}{k} \\
 &+& \left( \frac{3}{2} \, \frac{\Ut^1_h-\Ut^0_h}{k}-\frac{1}{2} \, \frac{\Ut^N_h-\Ut^{N-1}_h}{k} \, e^{-2i\pi T}
 \right) \frac{e^{-2i\pi k}-1}{k}.
\end{eqnarray*}
Taking the scalar product  with $i \, \widehat \Ut_{\eps}(\tau)$ we
get
\begin{equation*}
\label{eq:eqconv1505} 2   \pi  \tau \, \left( \frac{3}{2}
-\frac{1}{2} \, e^{-2 i\pi \tau k}\right) |\widehat
\Ut_{\eps}(\tau)|^2 =i\left( \mathbf{\Delta}_h \widehat
\G_{\eps}(\tau),\widehat \Ut_{\eps}(\tau) \right) +i
\left(\alphag,\widehat \Ut_{\eps}(\tau)\right).
\end{equation*}
Let us bound the right-hand side. According to  proposition \ref{prop:coerlap} and (\ref{eq:conveq135})
\begin{equation*}
\left|i\left( \mathbf{\Delta}_h \widehat \G_{\eps}(\tau),\widehat
\Ut_{\eps}(\tau) \right)\right| \le \|\widehat \G_{\eps}(\tau)\|_h
\, \|\widehat \Ut_{\eps}(\tau)\|_h
\le C \, \|\widehat \Ut_{\eps}(\tau)\|_h.
\end{equation*}
On the other hand, using theorem \ref{theo:stab}, one checks that $\alphag$ is bounded. Thus, according
to the Cauchy-Schwarz inequality and (\ref{eq:inpoinp0})
\begin{equation*}
\left|i \left(\alphag,\widehat \Ut_{\eps}(\tau)\right)\right| \le
|\alphag| \, |\widehat \Ut_{\eps}(\tau)| \le C \, |\alphag| \,
\|\widehat \Ut_{\eps}(\tau)\|_h\le C  \, \|\widehat
\Ut_{\eps}(\tau)\|_h .
\end{equation*}
Hence we have
\begin{equation*}
\label{eq:conveq18}
  \forall \, \tau \in \mathbb{R}, \, \hspace{1cm} |\tau| \, |\widehat \Ut_{\eps}(\tau)|^2 \le C \, \|\widehat \Ut_{\eps}(\tau)\|_h.
\end{equation*}
If $\tau \neq 0$, multiplying this estimate by
$|\tau|^{2 \, \gamma-1}$, we get
 $|\tau|^{2 \, \gamma} \, |\widehat \Ut_{\eps}(\tau)|^2 \le
 C \,
 |\tau|^{2 \, \gamma-1} \, \|\widehat \Ut_{\eps}(\tau)\|_h$.
Using the Young inequality and integrating over $\{\tau \in
\mathbb{R} \, ; |\tau| > 1\}$ we obtain
\begin{equation*}
\int_{|\tau| > 1} |\tau|^{2 \, \gamma} \, |\widehat
\Ut_{\eps}(\tau)|^2 \, d\tau \le \int_{|\tau| > 1}  |\tau|^{4 \,
\gamma-2} \, d\tau + C \, \int_{|\tau| > 1} \|\widehat
\Ut_{\eps}(\tau)\|^2_h \, d\tau.
\end{equation*}
For $|\tau| \le 1$ we have $ |\tau|^{2 \, \gamma} \, |\widehat
\Ut_{\eps}(\tau)|^2 \le |\widehat \Ut_{\eps}(\tau)|^2 \le C \,
\|\widehat \Ut_{\eps}(\tau)\|^2_h $ thanks to (\ref{eq:inpoinp0}).
Thus
\begin{equation*}
\int_{\mathbb{R}} |\tau|^{2 \, \gamma} \, |\widehat
\Ut_{\eps}(\tau)|^2 \, d\tau \le \int_{|\tau| > 1}  |\tau|^{4 \,
\gamma-2} \, d\tau+ C \, \int_{\mathbb{R}} \|\widehat
\Ut_{\eps}(\tau)\|^2_h \, d\tau.
\end{equation*}
Since $4 \, \gamma-2<-1$ we have $\int_{|\tau| > 1}  |\tau|^{4 \,
\gamma-2} \, d\tau \le C$. On the other hand, thanks to the Parseval
theorem and thorem \ref{theo:stab}
\begin{equation*}
  \int_{\mathbb{R}} \|\widehat \Ut_{\eps}(\tau)\|^2_h \, d\tau \le  \int_{\mathbb{R}} \|\widehat \Ut_{\eps}(\tau)\|^2_h \, dt
 \le k \sum_{n=1}^N \|\Ut^n_h\|^2_h \le C.
\end{equation*}
Hence the result.
\qed


\noindent We introduce the following spaces
\begin{equation*}
  \mathbf{H}=\{ \V \in \L^2 \, ; \hspace{.1cm} \hbox{div} \, \V  \in L^2 \hspace{.1cm}
  \hbox{ et } \hspace{.1cm} \V \cdot \N|_{\partial \Omega}=0\}  \, ,
  \hspace{.6cm}
  \mathbf{V}=\{ \V \in \H^1_0 \, ; \hspace{.1cm} \hbox{div} \, \V=0\}.
\end{equation*}
We also set
\begin{equation*}
((\U,\V))=\sum_{i=1}^2 (\nabla u_i,\nabla v_i) \, , \hspace{1cm}
\b(\U,\V,\W)=-\sum_{i=1}^2 \left(v_i, \U\cdot \nabla w_i\right)
\end{equation*}
for
all $\U=(u_1,u_2) \in \H^1$, $\V=(v_1,v_2) \in \H^1$, $\W=(w_1,w_2) \in \H^1$.
 We have the
following result.
\begin{theorem}
\label{theo:conv}

\noindent We assume that the initial values of the scheme fulfill
hypothesis {\bf (HI)}. We also assume that the space step $h$ and
the time step $k$ are such that
 $h \le C \, k^\alpha$ with $\alpha>1$.
 Then  we have $\U_\eps \to \U$ in $L^2(0,T;\L^2)$ with
\begin{equation}
\label{eq:regu}
 \U \in {\cal{C}}(0,T;\mathbf{H}) \cap L^2(0,T;\mathbf{V}) \, ,
 \hspace{1.5cm} \frac{d\U}{dt} \in L^2(0,T;\L^2).
\end{equation}
We also have $\U(0)=\U_0$ and for all $\psi \in
{\cal{C}}^{\infty}_0([0,T])$
\begin{equation}
\label{eq:formvarconv}
 \forall \, \V \in \mathbf{V}, \hspace{.5cm}
 \int^T_0    \psi  \left(\frac{d}{dt}(\U,\V) +
((\U,\V))+\b(\U,\U,\V)-(\F,\V) \right)  dt=0.
\end{equation}
\end{theorem}
\noindent {\sc Proof.} In what follows, sub-sequences  of a
sequence $(\V_\eps)_{\eps>0}$ will still be noted
$(\V_\eps)_{\eps>0}$ for the sake of convenience. All the limits are for  $\eps \to 0$.
According to theorem \ref{theo:stab} and hypothesis {\bf (HI)} we
have
\begin{equation*}
  \|\U_{\eps}\|^2_{L^2(0,T;\L^2)}= k \, (|\U^0_h|^2+|\U^1_h|^2)+k\sum_{n=2}^N |\U^n_h|^2 \le C.
\end{equation*}
We also deduce from (\ref{eq:inpoinp0}), hypothesis {\bf (HI)} and
theorem \ref{theo:stab}
\begin{equation*}
\|\Ut_{\eps}\|^2_{L^2(0,T;\L^2)} = k \, |\U^1_h|^2+k\sum_{n=2}^N
|\Ut^n_h|^2
  \le C+C \, k\sum_{n=2}^N \|\Ut^n_h\|^2_h \le C.
\end{equation*}
A simple computation shows that there exists $C>0$ such that
\begin{equation*}
\|\Ut^c_{\eps}\|_{L^2(0,T;\L^2)} \le C \,
\|\Ut_{\eps}\|_{L^2(0,T;\L^2)} \le C.
\end{equation*}
Thus the sequences  $(\U_{\eps})_{\eps>0}$, $(\Ut_{\eps})_{\eps>0}$
and $(\Ut^c_{\eps})_{\eps>0}$ are bounded
 in $L^2(0,T;\L^2)$.
 Therefore there exists  $\U \in L^2(0,T;\L^2)$, $\Ut \in L^2(0,T;\L^2)$
 and
$\Ut^{c} \in L^2(0,T;\L^2)$ such that, up to a sub-sequence, we have
\begin{equation*}
\label{eq:conveq21} \U_{\eps} \rightharpoonup  \U \, , \hspace{.5cm}
\Ut_{\eps} \rightharpoonup  \Ut \, ,  \hspace{.5cm} \Ut^c_{\eps}
\rightharpoonup  \Ut^{c} \hspace{.5cm} \hbox{ weakly in  }L^2(0,T;
\L^2).
\end{equation*}
We claim that the limits  $\U$, $\Ut$, $\Ut^{c}$ are the same. Indeed, let us
consider $\U_{\eps} - \Ut_{\eps}$. Since
$\U^n_h - \Ut^{n+1}_h=(\U^n_h - \U^{n+1}_h)+(\U^{n+1}_h -
\Ut^{n+1}_h)$
we have
\begin{equation*}
\label{eq:conveq25}
\|\U_{\eps} - \Ut_{\eps}\|^2_{L^2(0,T;\L^2)} 
\le 2 \, k\sum_{n=0}^{N-1} |\U^n_h - \U^{n+1}_h|^2+2 \,
k\sum_{n=0}^{N-1} |\U^{n+1}_h - \Ut^{n+1}_h|^2.
\end{equation*}
According to  theorem \ref{theo:stab} we have
$k\sum_{n=0}^{N-1} |\U^n_h - \U^{n+1}_h|^2 \le  C \sum_{n=0}^{N-1} k^3 \le C \, k^2$.
Thanks to (\ref{eq:projib}) we also have
\begin{equation*}
k\sum_{n=0}^{N-1} |\U^{n+1}_h - \Ut^{n+1}_h|^2= \frac{4 }{9} \,
k^3\sum_{n=1}^{N-1} |\nabla_h(p^{n+1}_h - p^n_h)|^2 \le C
\sum_{n=1}^{N-1} k^3 \le C \, k^2.
\end{equation*}
Thus $ \label{eq:conveq28}
  \|\U_{\eps} - \Ut_{\eps}\|_{L^2(0,T;\L^2)} \to 0
$ and
$ \label{eq:conveq282}
  \U=\Ut
$.
One checks in a simililar  way that
 $\label{eq:conveq289}
  \Ut=\Ut^c
  $.
 Now, using the Fourier transform, we prove the
 strong convergence  of the sequence $(\U_\eps)_{\eps>0}$ in $L^2(0,T;\L^2)$.
  We set
$ \V_\eps=\U_\eps-\U $.
  Let $M>0$. We use the splitting
\begin{equation*}
\label{eq:conveq2023}
  \int_{\mathbb{R}} |\widehat{\V}_\eps(\tau)|^2 \, d\tau
  = \int_{|\tau| \le M} |\widehat{\V}_\eps(\tau)|^2 \, d\tau
  + \int_{|\tau| > M} |\widehat{\V}_\eps(\tau)|^2 \, d\tau
=I^M_{\eps}+J^M_{\eps}.
\end{equation*}
Let us estimate $J^M_\eps$. Since
\label{eq:conveq202} $  |\widehat{\V}_\eps(\tau)|^2
   \le 2 \, |\widehat{\U}_\eps(\tau)|^2 +2 \,
  |\widehat{\U}(\tau)|^2$
we have
\begin{equation*}
\label{eq:conveq203} J^M_{\eps} \le 2 \, \int_{|\tau| > M}
|\widehat{\U}_\eps(\tau)|^2 \, d\tau + 2 \, \int_{|\tau| > M}
|\widehat{\U}(\tau)|^2 \, d\tau.
\end{equation*}
According to  lemma \ref{lemcompvi} we have
\begin{equation*}
\label{eq:conveq204} \int_{|\tau|> M} |\widehat{\U}_\eps(\tau)|^2 \,
d\tau \le \frac{1}{M^{2\gamma}} \int_{|\tau| > M}
  |\tau|^{2\gamma} \, |\widehat{\U}_\eps(\tau)|^2 \, d\tau
\le \frac{C}{M^{2 \, \gamma}}.
\end{equation*}
Thus
\begin{equation*}
\label{eq:conveq2055} J^M_{\eps} \le  \frac{2 \, C}{M^{2 \, \gamma}}
+2\int_{|\tau| > M} |\widehat{\U}(\tau)|^2 \, d\tau.
\end{equation*}
Therefore, for all $\eps>0$, we have $J^M_\eps \to 0$ when $M \to
\infty$. We now consider $I^M_\eps$.
Let $\tau \in \mathbb{R}$. Since $\U_\eps \rightharpoonup \U$ in $L^2(0,T;\L^2)$, we deduce from
 definition (\ref{eq:deftf})
\begin{equation*}
\label{eq:conveq2065}
  \widehat{\Ut}_\eps(\tau) \rightharpoonup \widehat{\U}(\tau) \hspace{.5cm}\hbox{ weakly in  $\L^2$}.
\end{equation*}
For all $t \in \mathbb{R}$ we have  $\Ut_\eps(t) \in \P_0$. From
definition (\ref{eq:deftf}) we infer that  $\widehat{\Ut}_\eps(\tau)
\in \P_0$. Now, prolonging $\widehat{\Ut}_\eps(\tau)$ by $0$ outside $\Omega$, we deduce
 from lemma 4 in \cite{eymgal} that  there exists a
constant $C>0$ such that
\begin{equation*}
\label{eq:conveq2066}
\forall \, \Eta{}{} \in \mathbb{R}^2\, , \hspace{1cm}   |\widehat{\Ut}_\eps(\tau)(\cdot + \Eta{}{})
- \widehat{\Ut}_\eps(\tau)|^2 \le
\|\widehat{\Ut}_\eps(\tau)\|^2_h \, |\Eta{}{}| \, (|\Eta{}{}| +C \,
h).
\end{equation*}
Using definition (\ref{eq:deftf}), the Cauchy-Schwarz inequality
 and theorem \ref{theo:stab},  we have
\begin{equation*}
\label{eq:conveq207}
 \|\widehat{\Ut}_\eps(\tau)\|^2_h \le C \int^T_0 \|\Ut_\eps(t)\|^2_h \, dt \le
 C \, k \sum_{n=1}^N \|\Ut^n_h\|^2_h \le C.
\end{equation*}
Thus, using the compactness criterium given by theorem 1 in \cite{eymgal}, we get
$ \label{eq:conveq20801} \widehat{\Ut}_\eps(\tau) \to \widehat{\U}(\tau)$ in
$\L^2$. Thus $ \widehat{\widetilde{\V}}_\eps(\tau)= \widehat{\Ut}_\eps(\tau)-
\widehat{\U}(\tau) \to 0$ in $\L^2$. Therefore for all $M>0$ we have
$I^M_\eps \to 0$.
Using the Parseval inequality, and gathering the limits for
$I^M_\eps$ and $J^M_\eps$, we get
\begin{equation*}
\label{eq:conveq2085} \int_{\mathbb{R}} |\widehat{ \V}_\eps(\tau)|^2
\, d\tau = \int_{\mathbb{R}} |{\V}_\eps|^2 \, dt= \int_{\mathbb{R}}
|{\U}_\eps-\U|^2
 \, d\tau \to 0.
\end{equation*}
We have proven that $\U_\eps \to \U$ in $L^2(0,T;\L^2)$. 

\noindent We now check the properties of $\U$.
First, proceeding as in  \cite{eymgal}, one checks  easily that $ \U \in
L^2(0,T;\H^1_0)$.
Now let $q \in L^2(0,T;{\cal{C}}^\infty_0)$.
According to (\ref{eq:esterrintp1nc}) we have
  $\nabla_h (\Pi_{P^{nc}_1} q) \to \nabla q $ in $L^2(0,T;\L^2)$.
Since $\U_\eps \to \U$ in $L^2(0,T;\L^2)$ we get
\begin{equation*}
 \big( \nabla_h (\Pi_{P^{nc}_1} q),\U_\eps \big) \to (\nabla q,\U)=-(q,\hbox{div} \, \U).
\end{equation*}
On the other hand, according to
propositions \ref{prop:umrt0} and \ref{prop:propadjh}, we have for
all $\eps>0$
\begin{equation*}
  \big(\nabla_h(\Pi_{P^{nc}_1} q),\U_\eps \big)=-(\Pi_{P^{nc}_1}q,\hbox{div}_h \, \U_\eps)=0.
\end{equation*}
Thus we have
  $\int^T_0 q \,\hbox{div} \, \U \, dt=0$ for all   $q \in L^2(0,T;{\cal{C}}^\infty_0)$.
Since the space ${\cal{C}}^\infty_0$ is dense in  $L^2$, we get
$
\label{eq:conveq489}
  \hbox{div} \, \U =0
$. Hence $ \label{eq:conveq4899}
 \U \in L^2(0,T;\mathbf{V})
$. \noindent Let us now check the  regularity of $\displaystyle
\frac{d\U}{dt}$. Using hypothesis {\bf (HI)}, (\ref{eq:inpoinp0}) and
theorem \ref{theo:stab}, we have
\begin{equation*}
  \left\|\frac{d\Ut^c_\eps}{dt}\right\|^2_{L^2(0,T;\L^2)}= \frac{1}{k} \, |\U^1_h-\U^0_h|^2+k \sum_{n=2}^N
  \frac{|\Del \Ut^n_h|^2}{k^2} \le C+C \, k \sum_{n=2}^N  \frac{\|\Del \Ut^n_h\|^2_h}{k^2}
  \le C.
\end{equation*}
Thus the sequence  $\displaystyle \left(
\frac{d\Ut^c_\eps}{dt}\right)_{\eps>0}$ is bounded in
$L^2(0,T;\L^2)$. Since $\U_\eps \to \U$ in $L^2(0,T;\L^2)$ with $\U \in L^2(0,T;\H^1)$, proceeding as
in, we get
\begin{equation*}
\frac{d\Ut^c_\eps}{dt} \rightharpoonup \frac{d\U}{dt} \hbox{ weakly }\in L^2(0,T;\L^2)
\end{equation*}
  and $\U \in {\mathcal{C}}(0,T;\H)$.

\noindent Let us now prove that $\U$ satisfies
(\ref{eq:formvarconv}). For the sake of simplicity, we omit to note
some time dependencies. According to (\ref{eq:mombdf}) we have for
all $t \in [t_1,T]$
\begin{eqnarray*}
\label{eq:conveq333}  &&  \frac{3}{2} \, \frac{d\Ut^c_{\eps}}{dt}
-\frac{1}{2} \, \frac{d\Ut^c_{\eps}}{dt}(t-k)
 - \frac{1}{\hbox{Re}} \, \mathbf{\Delta}_h \Ut_{\eps} +\bt_h\big(2 \, \U_\eps-\U_\eps(t-k),\Ut_\eps\big)
  -\F_\eps \\
&&= - \left( 1+\frac{4}{3} \, \chi_{[t_3,T]}\right) \nabla_h p_\eps
+\frac{5}{3} \, \chi_{[t_3,T]} \, \nabla_h p_\eps(t-k)-\frac{1}{3}
\, \chi_{[t_3,T]} \, \nabla_h p_\eps(t-2 \, k).
\end{eqnarray*}
Let $\V \in \mathbf{V} \cap ({\cal{C}}^\infty_0)^2$ and $\psi \in
{\cal{C}}^\infty([0,T])$ with $\psi(T)=0$. We set
$\V_h=\widetilde \Pi_{\P_0} \V$. Multiplying the former equation
by $\psi \, \V_h$ and integrating over $[t_1,T]$
 we get
\begin{eqnarray}
\label{eq:conveq34} &&  \int^T_{t_1}  \psi   \left( \frac{3}{2} \,
\frac{d\Ut^c_\eps}{dt} -\frac{1}{2} \frac{d\Ut^c_\eps}{dt}(t-k) ,
\V_h\right)  dt
 - \frac{1}{\hbox{Re}}\int^T_{t_1} \psi\, (\Deltat_h \Ut_\eps,\V_h)  \, dt \nonumber \\
&&
 +\int^T_{t_1} \psi \, \b_h\big(2 \, \U_\eps-\U_\eps(t-k),\Ut_\eps,\V_h\big)  \, dt
-\int^T_{t_1} \psi \, (\F_\eps,\V_h) \, dt \nonumber \\
&& = \int^T_{t_1} \chi \, \psi \, ( \nabla_h p_\eps,\V_h) \, dt
\end{eqnarray}
with
  $\chi=-\chi_{[t_3,T-2 \, k]}+\frac{1}{3} \,\chi_{[t_2,t_3]}-\frac{4}{3} \, \chi_{[t_1,t_2]}
  -\frac{7}{3} \, \chi_{[T-k,T]}-\frac{2}{3} \, \chi_{[T-2 \, k,T-k]}$.
We now check the limits of the terms in this equation.
First, according to (\ref{eq:esterrintp0}), we have $\V_h \to \V$ in $\L^2$. We will use this limit in
 the computations below without  mentioning it.
Since $\psi(T)=0$ we obtain by integrating by parts
\begin{equation*}
\label{eq:conveq49}   \int^T_{t_1}  \psi
\left(\frac{d\Ut^c_{\eps}}{dt}, \V_h\right)  dt-\psi(t_1) \, (\Ut^1_{h},\V_h)-\int^T_{t_1} \psi' \,
(\Ut^c_{\eps},\V_h) \, dt
\end{equation*}
and
\begin{equation*}
 \int^T_{t_1}  \psi \left(\frac{d\Ut^c_{\eps}}{dt}(t-k), \V_h\right)  dt
 =-\psi(0) \, (\Ut^0_{h},\V_h)-\int^{T-k}_0 \psi'(t+k) \, (\Ut^c_{\eps},\V_h) \, dt.
\end{equation*}
According to hypothesis {\bf (HI)} we have
  $\Ut^0_h=\U^0_h \to \U_0$ in $\L^2$ and  $\Ut^1_h=\U^1_h \to \U_0$ in $\L^2$. It implies that
  $(\U^0_h,\V_h) \to (\U_0,\V)$ and
   $\psi(t_1) \, (\Ut^1_{h},\V_h)= \psi(k) \, (\Ut^1_{h},\V_h)\to \psi(0) \, (\U_0,\V)$.
On the other hand
\begin{equation*}
  \int^T_{t_1} \psi' \, (\Ut^c_{\eps},\V_h) \, dt
  = \int^T_{0} \chi_{[t_1,T]} \, \psi' \, (\Ut^c_{\eps},\V_h) \, dt
  \to \int^T_0 \psi' \, (\U, \V) \, dt
\end{equation*}
and since $\chi_{[0,T-k]} \, \psi'(\cdot+k) \to \psi'$ in
$L^\infty(0,T)$
\begin{equation*}
\int^{T-k}_{0} \! \!  \psi'(t+k) \, (\Ut^c_{\eps},\V_h) \, dt
  = \int^T_0 \! \! \chi_{[0,T-k]} \, \psi'(t+k) \, (\Ut^c_{\eps},\V_h) \, dt
  \to \int^T_0 \! \psi' \, (\U, \V) \, dt.
\end{equation*}
Thus we have
\begin{equation}
\label{eq:conveq511}
  \int^T_{t_1}  \! \psi \left( \frac{3}{2} \frac{d\Ut^c_{\eps}}{dt}- \frac{1}{2} \, \frac{d\Ut^c_{\eps}}{dt}(t-k),
  \V_h\right)  dt \to
    -\psi(0) \, (\U,\V)-\int^T_0 \! \psi' \, (\U,\V) \, dt.
\end{equation}
Let us now consider the discrete laplacian. Using
proposition \ref{prop:adjlap} and
 the  splitting
 $\Deltat_h \V_h=
\big(\Deltat_h \V_h - \Pi_{\P_0} (\Deltat \V)\big) +
\Pi_{\P_0} (\Deltat \V)$
we have
\begin{equation*}
\label{eq:conveq355} \int^T_0  \psi \,  \big(\Deltat_h
\Ut_{\eps},\V_h)  \, dt
  = A_\eps+B_\eps
\end{equation*}
with
$A_\eps=\int^T_0 \psi \, \big( \Ut_{\eps}, \Deltat_h ( \widetilde
\Pi_{\P_0} \V)- \Pi_{\P_0}(\Deltat \V) \big)  \, dt$,
 $B_\eps=\int^T_0 \psi \, \big( \Ut_{\eps}, \Pi_{\P_0}(\Delta \V)
\big)  \, dt$.
Since
\begin{equation*}
\label{eq:conveq37} |A_\eps| \le  \|\Deltat_h (\widetilde \Pi_{\P_0}
\V )-\Pi_{\P_0}(\Deltat \V)\|_{-1,h} \int^T_0 \psi \, \|\Ut_\eps\|_h
\, dt \, ,
\end{equation*}
using proposition \ref{propconslap} and the Cauchy-Schwarz
inequality, we get
\begin{equation*}
|A_\eps| \le C \, h \, \|\V\|_2 \, \int^T_0 \psi \, \|\Ut_\eps\|_h
\, dt\le C \, h  \, \left( \int^T_0 \psi^2 \, dt \right)^{1/2}
  \left( \int^T_0 \|\Ut_\eps\|^2_h \, dt\right)^{1/2}.
\end{equation*}
Therefore, using theorem \ref{theo:stab}:
 $|A_\eps|
  \le C \, h  \left( k \sum_{n=1}^N \|\Ut^n_h\|^2_h\right)^{1/2} \le C \, h$.
Hence
$
\label{eq:conveq371}
  A_\eps \to 0
$. On the other hand, using an integration by parts, we have
\begin{equation*}
\label{eq:conveq383} B_\eps=\int^T_0 \psi \, ( \Ut_{\eps}, \Delta \V
)  \, dt \to \int^T_0 \psi \, (\U,\Delta \V) \, dt=
- \int^T_0 \psi \, ((\U,\V)) \, dt.
\end{equation*}
By gathering the limits for $A_\eps$ and $B_\eps$ we get
\begin{equation*}
\label{eq:conveq70}
  \int^T_0 \psi \, (\Deltat_h \Ut_{\eps},\V_h) \, dt \to
 - \int^T_0 \psi \, ((\U,\V)) \, dt.
\end{equation*}
Let us now consider the pressure. We use the splitting 
\begin{equation}
\label{eq:decompgh} (\nabla_h p_\eps,\V_h)= ( \nabla_h
p_\eps,\V_h-\V )
  +(\nabla_h p_\eps,\V-\Pi_{\mathbf{RT_0}} \V)
  +(\nabla_h p_\eps,\Pi_{\mathbf{RT_0}} \V).
\end{equation}
First, integrating by parts, we have
\begin{equation*}
(\nabla_h p_\eps,\Pi_{\mathbf{RT_0}} \V) =-\big( p_\eps, \hbox{div}
\, (\Pi_{\mathbf{RT_0}} \V) \big) +\sum_{K\in{\cal{T}}_h}
\int_{\partial K} p_\eps \, (\Pi_{\mathbf{RT_0}} \V \cdot
\N_{K,\sigma}).
\end{equation*}
Since $\hbox{div} \, \V=0$, using the divergence formula and definition (\ref{eq:defprt0}),
 one checks that $\hbox{div} \, (\Pi_{\mathbf{RT_0}} \V)=0$.  Thus $-\big( p_\eps,
\hbox{div} \, (\Pi_{\mathbf{RT_0}} \V) \big)=0$.
 On the other hand
\begin{equation*}
\sum_{K\in{\cal{T}}_h} \int_{\partial K} p_\eps \,
(\Pi_{\mathbf{RT_0}} \V \cdot \N_{K,\sigma}) =\sum_{\sigma \in
{\cal{E}}^{int}_h} \big((\Pi_{\mathbf{RT_0}} \V)_\sigma \cdot
\N_{K_\sigma,\sigma})
 \int_\sigma (p_\eps|_{L_\sigma}-p_\eps|_{K_\sigma}) \, d\sigma
\end{equation*}
and since $p_\eps\in P^{nc}_1$ we get
$\sum_{K\in{\cal{T}}_h} \int_{\partial K} p_\eps \,
(\Pi_{\mathbf{RT_0}} \V \cdot \N_{K,\sigma}) =0$.
Thus the last term in (\ref{eq:decompgh}) vanishes.
To bound the other terms, we use the Cauchy-Schwartz inequality together with estimates (\ref{eq:esterrintp0}),
 (\ref{eq:esterrirt0})
and theorem \ref{theo:stab}.
We get
\begin{equation*}
\left|(\nabla_h p_\eps,\V-\V_h) \right| +\left|(\nabla_h
p_\eps,\V-\Pi_{\mathbf{RT_0}}\V) \right|
\le C \, h \, |\nabla_h p_\eps| \, \|\V\|_2
\le C \, \frac{h}{k} \,
\|\V\|_2.
\end{equation*}
Plugging these estimates into (\ref{eq:decompgh}) we get
   $\int^T_{t_1} \left| \chi \, \psi \, (\nabla_h p_{\eps},\V_h) \right| \, dt
  \le C \, \frac{h}{k}$.
By hypothesis we have $\frac{h}{k} \le k^{\alpha-1}$ with
$\alpha-1>0$. Thus for $\eps=\max(h,k) \to 0$
\begin{equation*}
\label{eq:conveq60}
  \int^T_{t_1} \chi \, \psi \, (\nabla_h p_{\eps}, \V_h) \, dt \to 0.
\end{equation*}
Let us now consider the convection term. We set
  $\overline \U_\eps=2 \, \U_\eps -\U_\eps(t-k)$
and want to find the limit of
$\int^T_{t_1} \psi \, \b_h(\overline \U_\eps,\Ut_\eps,\V_h) \, dt$.
We use the splitting
$- \b_h(\overline \U_\eps,\Ut_\eps, \V_h)
+\b(\U,\U,\V)=A^\eps_1+ A^\eps_2+A^\eps_3$
with
\begin{equation*}
A^\eps_1=\b(\U-\overline \U_\eps,\U,\V) \, , \hspace{1cm}
A^\eps_2= \b(\overline \U_\eps,\U,\V)-\sum_{i=1}^2 \big(\hbox{div}( u_i\, \overline \U_\eps) , v^i_h \big) \, ,
\end{equation*}
\begin{equation*}
A^\eps_3= \sum_{i=1}^2 \big(\hbox{div}( u_i \, \overline \U_\eps)
, v^i_h \big)- \b_h(\overline \U_\eps,\U_\eps,\V_h).
\end{equation*}
By definition
  $A^\eps_1  =-\sum_{i=1}^2 \big(u_i , (\U-\overline \U_\eps) \cdot  \nabla v_i \big)$.
Using the Cauchy-Schwarz inequality we get
\begin{equation*}
\label{eq:formvar40}
 \int^T_{t_1} \psi \, |A^\eps_1| \, dt \le \|\psi\|_\infty \, \|\V\|_{\mathbf{W}^{1,\infty}} \, \|\U\|_{L^2(0,T;\L^2)} \,
 \|\U-\overline \U_\eps\|_{L^2(0,T;\L^2)}.
\end{equation*}
Since $\U_\eps \to \U$ in $L^2(0,T;\L^2)$ we also have $\|\U-\overline \U_\eps\|_{L^2(0,T;\L^2)} \to 0$. Thus
  $\int^T_{t_1} \psi \, A^\eps_1 \, dt \to 0$.
Let us now consider $A^\eps_2$. Since $\overline \U_\eps \cdot
\N|_{\partial \Omega}=0$ we obtain by integrating by parts
  $\b(\overline \U_\eps,\U,\V)
  =\sum_{i=1}^2\big(v_i,\hbox{div}(u_i \, \overline \U_\eps)\big)$.
Thus
\begin{equation*}
  A^\eps_2=
\sum_{i=1}^2\big(v_i-v^i_h,\hbox{div}(u_i \, \overline \U_\eps)\big)
  =\sum_{i=1}^2\big(v_i-v^i_h,\overline \U_\eps \cdot \nabla u_i\big)
\end{equation*}
Using the Cauchy-Schwarz inequality we get
\begin{equation*}
\label{eq:formvar405}
  \int^T_{t_1}  \psi \, |A^\eps_2| \, dt \le C \, \|\psi\|_{L^\infty}\, \|\V-\V_h\|_{\L^\infty} \,
  \|\U_\eps\|_{L^2(0,T;\L^2)} \, \|\U\|_{L^2(0,T;\H^1)}.
\end{equation*}
Using a Taylor expansion one checks that
 $ \|\V-\V_h\|_{\L^\infty}  \le \|\V\|_{\mathbf{W}^{1,\infty}} \, h$.
We recall also that $\|\U_\eps\|_{L^2(0,T;\L^2)} \le C$.
Therefore
 $\int^T_{t_1} \psi \, A^\eps_2 \, dt \to 0$.
Let us now  bound $A^\eps_3$. For all triangle $K \in {\cal{T}}_h$
and all edge
 $\sigma \in {\cal{E}}_K \cap {\cal{E}}^{int}_h$, we set
$$
\widetilde \U^\eps_{K,L_\sigma}= \left\{
\begin{array}{lr}
  \widetilde \U_\eps|_ K & \hbox{if } (\overline \U_\eps \cdot \N_{K,\sigma}) \ge 0  \\
  \widetilde \U_\eps|_{L_\sigma} &  \hbox{if } (\overline \U_\eps \cdot \N_{K,\sigma}) < 0
\end{array}
\right. .
$$
Using the divergence formula one  checks that
\begin{equation*}
  A^\eps_3=\sum_{K \in {\cal{T}}_h} \sum_{\sigma \in {\cal{E}}_K \cap {\cal{E}}^{int}_h}
  \V_K \cdot \int_\sigma (\U-\widetilde \U^\eps_{K,L_\sigma}) \, (\overline \U_\eps \cdot \N_{K,\sigma}) \, d\sigma.
\end{equation*}
By writing this sum as a sum on the edges we get
\begin{equation*}
  A^\eps_3= \sum_{\sigma \in {\cal{E}}^{int}_h}
  (\V_{K_\sigma}-\V_{L_\sigma}) \cdot \int_\sigma (\U-\widetilde \U^\eps_{K_\sigma,L_\sigma}) \, (\overline \U_\eps \cdot \N_{K_\sigma,\sigma})\, d\sigma.
\end{equation*}
Thus, using definition (\ref{eq:defpp1nc}) and a
quadrature formula
\begin{eqnarray*}
  A^\eps_3&=& \sum_{\sigma \in {\cal{E}}^{int}_h}
  (\V_{K_\sigma}-\V_{L_\sigma}) \int_\sigma (\Pi_{\P^{nc}_1}\U-\widetilde \U^\eps_{K_\sigma,L_\sigma}) \,
   (\overline \U_\eps \cdot \N_{K_\sigma,\sigma})\, d\sigma \\
&=& \sum_{\sigma \in {\cal{E}}^{int}_h}
  (\V_{K_\sigma}-\V_{L_\sigma}) \,  |\sigma| \left( (\Pi_{\P^{nc}_1}\U)(\x_\sigma)-\widetilde \U^\eps_{K_\sigma,L_\sigma}\right)  (\overline \U_\eps \cdot \N_{K_\sigma,\sigma})_\sigma.
\end{eqnarray*}
We have $|\sigma| \le h$ and, using a Taylor expansion, one checks that
 $|\V_{K_\sigma}-\V_{L_\sigma}| \le  h \, \|\V\|_{\mathbf{W}^{1,\infty}}$.
Thus, thanks to the Cauchy-Schwarz inequality, we get
\begin{equation*}
  |A^\eps_3| \le C \, h^2 \left(\sum_{\sigma \in {\cal{E}}^{int}_h} |\overline \U_\eps(\x_\sigma)|^2 \right)^{1/2}
\left(\sum_{\sigma \in {\cal{E}}^{int}_h}  |(\Pi_{\P^{nc}_1}
\U)(\x_\sigma) -\widetilde \U^\eps_{K_\sigma,L_\sigma}|^2
\right)^{1/2}.
\end{equation*}
Using (\ref{eq:propairet}) we get
\begin{eqnarray*}
  |A^\eps_3| &\le& C \left( \sum_{K \in {\cal{T}}_h} \frac{|K|}{3}\sum_{\sigma \in {\cal{E}}_K \cap {\cal{E}}^{int}_h}
   |\overline \U_\eps(\x_\sigma)|^2 \right)^{1/2}  \\
&\times & \left( \sum_{K \in {\cal{T}}_h} \frac{|K|}{3}\sum_{\sigma
\in {\cal{E}}_K \cap {\cal{E}}^{int}_h} |(\Pi_{\P^{nc}_1}
\U)(\x_\sigma) -\widetilde \U_\eps|_{K}|^2 \right)^{1/2}.
\end{eqnarray*}
Therefore, using a quadrature formula
\begin{equation*}
\label{eq:formvar42}
 |A^\eps_3| \le C \, |\overline \U_\eps| \, |\Pi_{\P^{nc}_1} \U -\U_\eps|
\le C \, |\U_\eps| \, |\Pi_{\P^{nc}_1} \U -\U_\eps|.
\end{equation*}
By writing
$\Pi_{\P^{nc}_1} \U -\U_\eps= ( \Pi_{\P^{nc}_1} \U -\U ) +
(\U-\U_\eps)$
and using (\ref{eq:esterrintp1nc}), we get
 $|A^\eps_3| \le C \, |\U_\eps| \, (h \, \|\U\|_1+|\U-\U_\eps|)$.
Thus, using the  Cauchy-Schwarz inequality,  we have
\begin{equation*}
 \int^T_{t_1} |\psi| \, |A^\eps_3| \, dt \le C \,  \|\U_\eps\|_{L^2(0,T;\L^2)} \, ( h \, \|\U\|_{L^2(0,T;\H^1)}+
 \|\U-\U_\eps\|_{L^2(0,T;\L^2)}).
\end{equation*}
Thus
  $\int^T_{t_1} \psi \, A^\eps_3 \, dt \to 0$.
By gathering the limits for $A^\eps_1$, $A^\eps_2$, $A^\eps_3$  we obtain
  $\int^T_{t_1} \psi \, \big( \b_h(\overline \U_\eps,\U_\eps, \V_h)-\b(\U,\U,\V) \big)\, dt \to 0$.
Since
  $\int^{t_1}_0 \b(\U,\U,\V) \, dt  \to 0$, 
 we get
\begin{equation*}
\int^T_{t_1} \psi \, \big( \b_h(\overline \U_\eps,\U_\eps, \V_h) ,
dt \to \int_0^T \psi \, \b(\U,\U,\V) \, dt.
\end{equation*}
Finally, since $\V_h \in \P_0$, we have: $ \label{eq:conveq61}
(\F_{\eps}, \V_h)=(\Pi_{\P_0} \F, \V_h)=(\F,\V_h) $. Therefore
\begin{equation*}
\label{eq:conveq80}
 \int^T_{t_1} \psi \, (\F_{\eps}, \V_h) \, dt
= \int^T_{0} \chi_{[t_1,T]} \, \psi \, (\F, \V_h) \, dt
 \to \int^T_{0} \psi \, (\F, \V)  \, dt.
\end{equation*}
We now gather the limits we have obtained  into
(\ref{eq:conveq34}). The space  $\mathbf{V} \cap
({\cal{C}}^\infty_0)^2$ is dense in $\mathbf{V}$. Hence we obtain
%
 for all $\V \in \mathbf{V}$ and
$\psi \in {\cal{C}}^\infty([0,T])$ with $\psi(T)=0$
\begin{equation*}
\label{eq:conveq82}
 -\psi(0) \, (\U_0,\V)
 -\int^T_0 \psi' \, (\U,\V) \, dt
+ \int^T_0 \psi \, \big(((\U, \V))+\b(\U,\U,\V)-(\F,\V) \big) \,
dt=0.
\end{equation*}
Taking $\psi=\phi \in {\cal{C}}^{\infty}_0([0,T])$, we have
$\phi(0)=0$ and from the definition of the derivative in the
distributional sense
\begin{equation*}
\label{eq:conveq5101} \int^T_0 \phi' \, (\U,\V) \, dt =-\int^T_0
\phi \, \frac{d}{dt} (\U,\V) \, dt.
\end{equation*}
Thus we have proven (\ref{eq:formvarconv}).
\noindent At last, let us show that the initial condition
 holds. We have proven before that
\begin{equation*}
\label{eq:conveq9878}
  \frac{d\Ut^c_\eps}{dt} \rightharpoonup \frac{d\U}{dt} \hspace{.5cm} \hbox{ weakly  in }L^2(0,T;\L^2).
\end{equation*}
Let $\V \in \mathbf{V} \cap ({\cal{C}}^\infty_0)^2$ and $\psi \in
{\cal{C}}^\infty([0,T])$ such that $\psi(T)=0$.
We have
\begin{equation*}
\label{eq:conveq98785}
  \int^T_0  \psi \left(\frac{3}{2} \, \frac{d\Ut^c_\eps}{dt} - \frac{1}{2} \,\frac{d\Ut^c_\eps}{dt}(t-k) ,\V_h \right) \, dt
  \to \int^T_0 \psi \left(\frac{d\U}{dt}, \V \right)  \, dt.
\end{equation*}
Integrating by parts the limit we get
\begin{equation*}
\label{eq:conveq987991}
  \int^T_0  \psi \left(\frac{3}{2} \, \frac{d\Ut^c_\eps}{dt}-\frac{1}{2} \, \frac{d\Ut^c_\eps}{dt}(t-k),\V_h \right) \, dt
  \to -\psi(0) \, \big(\U(0),\V\big) -\int^T_0 \psi' \, (\U, \V) \, dt.
\end{equation*}
By comparing  this limit  with (\ref{eq:conveq511}),
we get
 $\psi(0) \, (\U(0)-\U_0,\V)=0$
for all $\psi \in {\cal{C}}^\infty([0,T])$ with
$\psi(T)=0$.
Therefore
  $\U(0)=\U_0$.
At last, note that we have proven so far the convergence of
a sub-sequence of $(\U_\eps)_{\eps>0}$ towards $\U$.
 But  the application $\U$
  such that (\ref{eq:regu}), (\ref{eq:formvarconv})
  and $\U(0)=\U_0$ hold is unique (\cite{temam}, p. 254).
Thus the whole sequence
$(\U_\eps)_{\eps>0}$ converges towards $\U$. \qed


\end{document}